\documentclass[11pt]{article}

\usepackage{amssymb,psfig}

\newtheorem{theorem}{Theorem}[section]

\newtheorem{defi}{Definition}[section]
\newtheorem{lemma}{Lemma}[section]

\def\slfrac#1#2{\hbox{\kern.1em %
 \raise.5ex\hbox{\the\scriptfont0 #1}\kern-.11em %
 /\kern-.15em\lower.25ex\hbox{\the\scriptfont0 #2}}}

\newcommand{\eqn}[1]{(\ref{#1})}
\newcommand{\hsp}{\hspace*{\parindent}}

\newcommand{\eeq}{\end{equation}}
\newcommand{\beql}[1]{\begin{equation}\label{#1}}
\newcommand{\bsq}{{\vrule height .9ex width .8ex depth -.1ex }}

\newcommand{\af}{\alpha}

\newcommand{\Om}{\Omega}

\newcommand{\RR}{{\mathbb R}}
\newcommand{\bv}{{\bf v}}
\newcommand{\bB}{{\bf B}}
\newcommand{\bC}{{\bf C}}
\newcommand{\bc}{{\bf c}}
\newcommand{\by}{{\bf y}}
\newcommand{\bx}{{\bf x}}
\newcommand{\bw}{{\bf w}}

\newcommand{\bzero}{{\bf 0}}
\newcommand{\sP}{{\cal P}}
\newcommand{\sG}{{\cal G}}

\newcommand{\sC}{{\cal C}}

\newcommand{\sD}{{\cal D}}

\newcommand{\sV}{{\cal V}}

\makeatletter
\def\@sect#1#2#3#4#5#6[#7]#8{\ifnum #2>\c@secnumdepth
     \def\@svsec{}\else
     \refstepcounter{#1}\edef\@svsec{\csname the#1\endcsname.\hskip .75em }\fi
     \@tempskipa #5\relax
      \ifdim \@tempskipa>\z@
        \begingroup #6\relax
          \@hangfrom{\hskip #3\relax\@svsec}{\interlinepenalty \@M #8\par}%
        \endgroup
       \csname #1mark\endcsname{#7}\addcontentsline
         {toc}{#1}{\ifnum #2>\c@secnumdepth \else
                      \protect\numberline{\csname the#1\endcsname}\fi
                    #7}\else
        \def\@svsechd{#6\hskip #3\@svsec #8\csname #1mark\endcsname
                      {#7}\addcontentsline
                           {toc}{#1}{\ifnum #2>\c@secnumdepth \else
                             \protect\numberline{\csname the#1\endcsname}\fi
                       #7}}\fi
     \@xsect{#5}}
\makeatother

\catcode`\@=11
\renewcommand{\section}{
        \setcounter{equation}{0}
        \@startsection {section}{1}{\z@}{-3.5ex plus -1ex minus
        -.2ex}{2.3ex plus .2ex}{\large\bf}
        }
\catcode`@=12

\begin{document}
\begin{center}
{\Large {\bf  Bounds for Local Density of Sphere Packings \\and
the  Kepler Conjecture}} \\
\vspace{1.5\baselineskip}
{\em Jeffrey  C. Lagarias} \\
%\vspace*{.2\baselineskip}
%AT\&T Labs - Research \\
%Florham Park, NJ 07932-0971 \\
%jcl@research.att.com \\
\vspace*{1.5\baselineskip}
\end{center}
\vspace{1.5\baselineskip}
{\small
 \centerline{\bf Abstract}

This paper formalizes the local density inequality approach
to getting upper bounds for sphere packing densities in $\RR^n$.
This approach was first suggested by L. Fejes-T\'oth 
in 1954 as a method
to prove the Kepler conjecture that the densest packing of
unit spheres in $\RR^3$ has density  $\frac{\pi}{\sqrt{18}}$,
which is attained by the ``cannonball packing.'' 
Local density inequalities give upper bounds for the sphere
packing density formulated as an optimization
problem of a nonlinear function over a compact set in a finite
dimensional Euclidean space. The approaches of 
L. Fejes-T\'oth, of W.-Y. Hsiang,  and of T. C. Hales,
to the Kepler conjecture are each based on 
(different) local density inequalities. Recently 
T. C. Hales, together with S. P. Ferguson,  has 
presented extensive details carrying out a modified
version of the Hales
approach to prove the Kepler conjecture. We describe
the particular local density inequality underlying the
Hales and Ferguson approach to prove Kepler's
conjecture and sketch some features of their proof.\\


}

\section{Introduction}
\hsp
The Kepler conjecture was stated by Kepler in 1611.
It asserts that the
face-centered cubic lattice gives the
tightest possible packing of unit spheres in $\RR^3$.

\paragraph{Kepler Conjecture.}
{\em Any packing $\Omega$ of unit spheres in $\RR^3$ has upper 
packing density
\beql{101}
\bar{\rho} (\Omega ) \le \frac{\pi}{\sqrt{18}}
\simeq 0.740480 ~.
\eeq
}

\vspace*{+.1in}
\noindent The definition of upper packing density is given in \S2.
The problem of proving the Kepler conjecture appears as part of 
Hilbert's 18th problem, see \cite{Hi}.

T. C. Hales has described an approach for proving 
Kepler's conjecture, and has announced a proof,
completed with the aid of S. Ferguson, which is currently presented in
a set of six preprints. The proof is computer-intensive,
and involves checking over $5000$ subproblems.
The Hales approach is similar to earlier approaches in that
it aims to prove a local density inequality that gives a (sharp)
upper bound on the density. It involves several new ideas
which are indicated in \S4 and \S5.

Local density inequalities 
obtain upper bounds for the sphere packing constant via an
auxiliary nonlinear optimization problem over a compact set of
``local configurations''. They measure a ``local density'' in
the neighborhood of each sphere center separately.
The general approach to the Kepler conjecture
is first to find a local optimization problem that 
actually attains the optimal bound
$\frac{\pi}{\sqrt{18}}$ (assuming that one exists), and then to prove it.
This approach was first suggested in the early 1950's by
L.  Fejes-T\'oth\cite[pp. 174-181]{FT}, who presented some
evidence that an optimal local density  inequality might exist
in three dimensions.

The objects of this paper are: 

(i) To formulate local density inequalities for sphere
packings in $\RR^n$, in sufficient generality to
include the known candidates for optimal
local inequalities.

(ii) To review the history of local density inequalities for
three dimensional sphere packing and the  Kepler conjecture.

(iii) To give a precise statement of the local density
inequality considered in the Hales-Ferguson  approach. 

(iv) To outline some of features of the Hales-Ferguson proof.

In \S2 we present a general framework for local density
inequalities, which is valid in $\RR^n$, 
given as Theorem~\ref{th21}. This framework is
sufficient to cover the approaches of L. Fejes-T\'oth,
W.-Y. Hsiang, and T. Hales and S. Ferguson to Kepler's conjecture.
A different framework for local density inequalities
appears in  Oesterl\'{e}~\cite{Oes99}.
In \S3 we review the history of work on local optimization inequalities
for Kepler's conjecture.
In \S4 we describe the precise local optimization problem 
formulated by Ferguson and Hales in \cite{FH},
which putatively attains $\frac{\pi}{\sqrt{18}}$.
In \S5 we make remarks on some details of the proof strategy  
taken in the papers of  Hales
\cite{I}, \cite{II}, \cite{III}, \cite{IV},
Hales and  Ferguson \cite{FH}, and Ferguson \cite{Fe}. 
In \S6 we make some concluding remarks.

The current status of the Hales-Ferguson proof is that 
it appears to be sound. The proof has reputedly been examined in
fairly careful detail by a team of reviewers, but it is so long
and complicated that it seems difficult for any one person to
check it. This paper is intended as an aid in understanding the 
overall structure of the Hales-Ferguson proof approach. For another
account of the Hales and Ferguson work,
see Oesterl\`{e} \cite{Oes99}.
For Hales' own perspective, see Hales \cite{H00}.

Two appendices are included which contain some information relevant to the
Hales-Ferguson proof.
Appendix A describes some of the Hales-Ferguson scoring functions.
Appendix B lists references in the Hales and Ferguson preprints for proofs of 
lemmas and theorems stated without proof in \S4 and \S5.

This paper is a slightly revised version of the manuscript ~\cite{La99}. 

\paragraph{Notation.}
$\bB_n := B_n( \bzero;1) = \{ \bx \in \RR^n : \| \bx \| \le 1 \}$ is the 
unit n-sphere. It has volume 
$\kappa_n :=  \pi^{\frac{n}{2} }/\Gamma(\frac{n}{2} + 1)$,
with $\kappa_2 = \pi$ and $\kappa_3 = \frac {4\pi}{3}$. We let
$\bC_n (\bx, T) := \bx + [0,T]^n$ 
denote an n-cube of sidelength $T$, with sides parallel to the 
coordinate axes,
and  lowest corner at $\bx \in \RR^n$.

\section{Local Density Inequalities}
\hsp
In this section we present a general formulation of local density
inequalities.

We recall the standard definition
of sphere packing densities,
following Rogers \cite{Ro}.
Let $\Omega$ denote a set of 
unit sphere centers, so that
$\| \bv- \bv' \| \ge 2 $ for distinct
$\bv , \bv' \in \Omega.$

\begin{defi}\label{de1}
{\rm
(i) For a bounded region $S$ in $\RR^n$, and a sphere packing $\Omega + \bB_n$
specified by the sphere centers $\Omega$, 
the {\em density}
$\rho (S)= \rho(\Omega, S)$ of the packing in the region $S$ is
\beql{103}
\rho (S) := \frac{vol (S \cap (\Omega +\bB_n))}{vol (S)} ~.
\eeq
(ii) For $T> 0$ the {\em upper density} $\rho (\Omega ,T)$ is the 
maximum density of the packing $\Omega$ over all cubes of size $T$, i.e.
\beql{104}
\bar{\rho} (\Omega, T ) := \sup_{\bx \in \RR^3} \rho([0, T]^n + \bx ) ~.
\eeq
Then the {\em upper packing density} of $\Omega$ is
\beql{105}
\bar{\rho} (\Omega ) := \limsup_{T \to \infty} \bar{\rho} (\Omega , T) ~.
\eeq
(iii) The {\em sphere packing density} $\delta(\bB_n)$ of the ball $\bB_n$
of unit radius is
\beql{eq105a}
\delta(\bB_n) : =  \sup_{\Omega} \bar{\rho} (\Omega ).
\eeq
}
\end{defi}

\begin{defi}\label{de22}
{\rm
A sphere packing $\Omega$ is {\em saturated} if no new sphere centers 
can be added to it.}
\end{defi}

To obtain sphere packing bounds it  obviously suffices to 
study saturated sphere packings,
and in what follows we assume that all packings are saturated
unless otherwise stated.

\begin{defi}\label{de23}
{\rm
An {\em admissible partition rule}
is a rule assigning to each saturated packing $\Omega$ in $\RR^n$
a collection of closed sets
$\sP ( \Omega ) : = \{ R_\alpha = R_\af ( \Omega ) \}$ with the
following properties.

(i) {\em Partition}.
Each set $R_\af$ is a finite union of bounded convex polyhedra.
The sets $R_\af $ cover $\RR^3$ and have pairwise disjoint interiors.

(ii) {\em Locality}.
There is a positive constant $C$ (independent of $\Omega$) such that 
each region $R_\af$ has
\beql{M21}
diameter (R_\af) \le C ~.
\eeq
Each $R_\af$ is completely determined by the set of sphere centers
$\bw \in \Omega$ with
\beql{M22}
distance ( \bw , R_\af ) \le C ~.
\eeq
There are at most $C$ regions intersecting any cube of side 1.

(iii) {\em Translation-Invariance}.
The partition assigned to the translated packing $\Om' = \Om + \bx$ 
consists of the sets
$\{ R_\af ( \Omega ) + \bx \}$.
}
\end{defi}

\begin{defi}\label{de4}
{\rm An {\em admissible weight function} (or {\em admissible score function}) 
$\sigma$ for an admissible
partition rule in $\RR^n$ assigns to each region
$R_\af \in \sP ( \Omega )$ and each $\bv \in \Omega$ a real weight 
$\sigma (R_\af , \bv )$ which satisfies  $|\sigma (R_\af , \bv )| < C^*$
for an absolute constant $C^*$, and  which has
the following properties.

(i)
{\em Weighted Density Average.}
There are positive constants $A$ and $B$ (independent of $\Omega$)
such that
for each set $R_\af$,
\beql{M23}
\sum_{\bv \in \Om} \sigma ( R_\af, \bv ) = (A~ \rho ( R_\af ) - B)
vol ( R_\af ) ~,
\eeq
where
\beql{M24}
\rho (R_\af ) vol (R_\af ) = vol (R_\af \cap (\Om + \bB_n) )
\eeq
measures the volume covered in $R_\af$ by the sphere packing $\Om$
with unit spheres.

(ii) {\em Locality.}
There is an absolute constant $C$ (independent of $\Omega$) such that each
value $\sigma (R_\af , \bv )$ is completely determined by the set of 
sphere centers
$\bw \in \Omega$ with $\| \bw - \bv \| \le C$.
Furthermore
\beql{M25}
\sigma (R_\af, \bv ) = 0 \quad \mbox{if}\quad
dist (\bv , R_\af ) > C~.
\eeq

(iii) {\em Translation-Invariance.}
The weight function $\sigma'$ assigned to the translated packing
$\Omega' = \Omega + \bx$ satisfies
\beql{M26}
\sigma' (R_\af + \bx , \bv + \bx ) = \sigma (R_\af, \bv ) ~.
\eeq
}
\end{defi}

Note that this definition specifically allows negative weights.

The ``local density'' is measured by  the sum of the  weights associated 
to a given vertex $\bv$ in a saturated packing. 
\begin{defi}\label{Nde25}
{\rm (i) The {\em vertex $D$-star} (or {\em decomposition star}) $\sD (\bv )$ 
at a vertex $\bv \in \Omega$ consists of all sets
$R_\af \in \sP (\Omega )$ such that
$\sigma (R_\af , \bv ) \neq 0$.

(ii). The {\em total score} assigned to a vertex $D$-star $\sD (\bv )$ at 
$\bv \in \Omega$ is
\beql{M27}
Score( \sD (\bv ) ) :=
\sum_{R_\af \in \sD (\bv )} \sigma ( R_\af , \bv )~.
\eeq
}
\end{defi}

The total score at $\bv$ depends only on
regions entirely contained within distance $C$ of $\bv$.
Any admissible partition and weight function
$(\sP , \sigma )$ together yield a local inequality for the density of
sphere packings, as follows.
\begin{theorem}\label{th21}
Given an admissible partition in $\RR^n$ 
and weight function $(\sP , \sigma )$,
set
\beql{M27a}
\theta = \theta_{\sP, \sigma} (A,B) := \sup_{\Omega ~\mbox{saturated}}
\left ( \sup_{\bv \in \Omega}~~ Score( \sD (\bv ) ) \right) ~.
\eeq
and suppose that $\theta < \kappa_n A,$ where $\kappa_n$ is
the volume of the unit $n$-sphere.
Then the maximum sphere packing density satisfies 
\beql{M28}
 \delta(\bB) \le \frac{\kappa_n B}{\kappa_n A- \theta} ~.
\eeq
\end{theorem}
\paragraph{Remark.}(1) We let  $f(A,B, \theta ) := 
\frac{ \kappa_n B}{\kappa_n A- \theta}$  denote the packing density bound
as a function of the score constants $A$ and $B$. 
The sphere packing density bound actually  depends only on the
{\em score constant ratio} $\frac {B}{A}$, rather than 
on $B$ and $A$ separately, since 
$\theta$ is a homongeneous linear function of $A$ and $B$. This ratio
detemines the relative weighting of covered and uncovered volume used
in the inequality. 

(2) A natural approach to sphere packing bounds, used in
many previous upper bounds,  is to partition space
into pieces $R(\bv)$ corresponding to each sphere center $\bv$, 
with each piece containing the unit sphere around $\bv$, and
aims to establish an upper bound
\beql{M28a}
\rho (R(\bv)) =  \frac {\kappa_n}{vol (R(\bv))} \le \theta~,~~~~
all~~~ \bv \in \Omega.
\eeq
Then one obtains $\bar{\rho}(\Omega) \le \theta.$ An 
optimal sphere packing bound of this sort
must necessarily be {\em volume-independent,} in the sense that
if equality is to be attained at all local cells $R(\bv)$
simultaneously, then they must all have the same volume.
In contrast the inequality
of Theorem~\ref{th21} does take into account the volumes of the 
individual pieces in the vertex $D$-star, and this allows 
more flexibility in the
local density inequalities that can be constructed, which might
attain an optimal bound.

\paragraph{Proof of Theorem~\ref{th21}.}
We may assume that $\Omega$ is saturated. Given $T > 0$, and any
$\epsilon > 0$ 
we choose a point $\bx \in \RR^n$ which attains the density bound
$\bar{\rho} (\Om, T )$ on the cube $\bC_n (\bx, T)$ to within $\epsilon$. 
We evaluate the scores of all vertex $D$-stars
of vertices $\bv \in \Omega \cap \bC_n (\bx, T)$ in two ways.
First, by definition of $\theta$,
\beql{M209}
\sum_{\bv \in \Omega \cap \bC_n(\bx, T )} Score (\sD (\bv) ) \le \theta \# | \Omega \cap \bC_n (\bx , T ) | ~.
\eeq
However we also have
\begin{eqnarray}\label{M210}
\sum_{\bv \in \Omega \cap \bC_n(\bx, T)} Score (\sD (\bv ) ) & = &
\sum_{\bv \in \Omega \cap \bC_n (\bx , T)}
\left( \sum_\af \sigma (R_\af, \bv )\right) \nonumber \\
& = & \sum_{\af \atop R_\af \subseteq \bC_n (\bx, T)}
\left( \sum_{\bv \in \Om \cap \bC_n(\bx, T)} \sigma (R_\af , \bv) \right) +
O(T^{n - 1}) \nonumber \\
& = & \sum_{R_\af \subseteq \bC_n (\bx;T )}
(A~~ \rho ( R_\af ) - B)~ vol (R_\af ) + O(T^{n - 1} ) \nonumber \\
& = & \kappa_n A \# | \Omega \cap \bC_n (\bx,T) | - 
B~~vol ( \bC_n (\bx, T)) + O (T^{n - 1} ) \nonumber \\
& = & \kappa_n A \#| \Omega \cap \bC_n (\bx, T ) | -
BT^n + O(T^{n - 1} ) ~.
\end{eqnarray}
Here we use that fact that $\{R_\af\}$ partitions $\RR^n$, so covers the cube,
and the $O(T^{n - 1})$ error terms above occur because the counting is not
perfect within a constant distance $C$ of the boundary of the cube.
Combining these evaluations yields
$$ (\kappa_n A ~-~\theta)\#| \Omega \cap \bC_n (\bx, T ) | \leq ~B~T^n + 
O(T^{n - 1}.) $$
If $\theta < \kappa_n A$, then we can rewrite this as
\beql{M211}
\frac {\#| \Omega \cap \bC_n (\bx, T ) |}{T^n} \leq 
\frac {B}{\kappa_n A - \theta}~ + O(\frac {1}{T}).
\eeq
By assumption
\begin{eqnarray*}
\bar{\rho} ( \Om, T) - \epsilon & \leq  &
\frac{vol ( \bC ( \bx , T) \cap ( \Omega + \bB ))}{T^3} \\
& = & \kappa_n
\frac{\#| \Omega \cap \bC (\bx, T)|}{T^n} + O \left( \frac{1}{T} \right) ~.
\end{eqnarray*}
Together with \eqn{M211}, this yields
\beql{M212}
\bar{\rho} ( \Omega, T) - \epsilon \le \frac{\kappa_n B}{\kappa_n A- \theta}
+ O \left( \frac{1}{T} \right) ~,
\eeq
with an $O$-symbol constant independent of $\epsilon.$
Letting $\epsilon \to 0$ and then 
$T \to \infty$ gives the inequality for $\bar{\rho} (\Omega)$.
Since this holds for all saturated packings the result follows.~~~$\bsq$

Determining the quantity $\theta_{\sP,\sigma} (A,B)$ for fixed $A,B$ can be
viewed as a nonlinear optimization problem over a compact set.
The translation-invariance property of $(\sP , \sigma )$ allows
the supremum \eqn{M27a} to be taken over the smaller set with
 $\bv = \bzero$ and admissible $\Omega$ containing $\bzero$. 
The locality property shows that that the vertex D-star at ${\bf 0}$ is
completely determined by $\bw \in \Omega$ with $\| \bw \| \le C$ .
The set of such configurations of nearby sphere centers forms a 
compact set in the Euclidean topology.
Actually the  partition and weight functions may be discontinuous functions
of the locations of sphere centers, so the optimization problem above
is not genuinely over a compact set.
One must compactify the space of  
allowable vertex $D$-stars 
by allowing some sets of sphere centers to be assigned
more than one possible
vertex $D$-star. In practical cases there is a finite upper bound
on the number of possibilities.

\begin{defi}\label{de14}
{\em A local density inequality in $\RR^n$ is} optimal 
{\em if}
\beql{M213}
f(A,B, \theta_{\sP, \sigma} ) =
\delta(\bB_n) ~.
\eeq
\end{defi}

Optimal local density inequalities exist in one and two dimensions.
In discussing the three-dimensional
case,  we shall presume that $\delta(\bB_3) = \frac {\pi}{\sqrt{18}}$,
so that an optimal density inequality  in $\RR^3$
will refer to one achieving this value. 
The evidence indicates that there are many different possible optimal
local density inequalities in three dimensions, including
that of the Hales and Ferguson proof. 

There are currently four candidates for local density 
inequalities that may be optimal in three dimensions. The first is that of 
L. Fejes-Toth, described in \S3,
which uses averages over Voronoi domains, 
in which the score constant ratio $\frac {B}{A} = \frac {\pi}{\sqrt{18}}.$ 
The second is that of Hsiang \cite{Hs}, which is a modification
of the Fejes-T\'oth averaging, and uses the same score constant ratio.
The third is due to Hales \cite{I}, and is based on the Delaunay 
triangulation, using a modified
scoring rule described in \S3. The fourth is that given in 
Ferguson and Hales \cite{FH}, and
uses a combination of Voronoi-type domains and Delaunay simplices,
with a complicated scoring rule, described in \S4. 
In the latter two cases the score constant ratio is
$$\frac{B}{A} = \delta_{oct}= \frac{-3 \pi + 
12 \arccos ( \frac{1}{\sqrt{3}} )}{\sqrt{8}} \approx 0.720903.$$
In all four cases the compact
set of local configurations to be searched
has very high dimension.
Each sphere center has three degrees of freedom, and the number of 
sphere centers
involved in these methods to determine a vertex $D$-star seems
 to be  around 50, so the 
search space consists of components of dimension up to roughly 150.

It is unknown whether {\em optimal} local density inequalities exist for 
the sphere packing problem in $\RR^n$ in any dimension $n \ge 4$. 
In dimensions 4, 8 and 24 it seems plausible that the minimal volume
Voronoi cell in any sphere packing actually
occurs in the densest lattice packing. If so, the Voronoi cell  
decomposition would yield
an optimal local inequality in these dimensions, and the densest packing
would be a lattice packing in these dimensions. 
Another question asks: in which dimensions is the maximal 
sphere packing density
attained by a sphere packing whose centers
form a finite number of cosets of an $n$-dimensional lattice? 
Perhaps in such dimensions
an optimal local density inequality exists. The state of the
art in sphere 
packings in dimensions four and above is given in 
Conway and Sloane \cite{CS}.
\section{History}
\hsp
We survey results on local density inequalities in
three dimensions.
The work on  local density bounds was originally based on two
partitions of $\RR^3$ associated to a set $\Omega$ of sphere
centers: the Voronoi tesselation and the Delaunay triangulation.
Since they will play an important role, we recall their definitions.

\begin{defi}\label{de31}
{\em The} Voronoi domain  ({\em or} Voronoi cell) {\em of $\bv \in \Omega$ is}
\beql{301}
V_{vor}( \bv ) = V_{vor}( \bv , \Omega) := \{ \bx \in \RR^3 : \| \bx-\bv \| \le \| \bx-\bw \| \quad\mbox{for all}\quad
\bw \in \Omega \} ~.
\eeq
{\em The} Voronoi tesselation {\em for $\Omega$ is the set of Voronoi domains}
$\{V_{vor} (\bv) : \bv \in \Omega \}$.
\end{defi}

The Voronoi tesselation is a partition of space, up to boundaries of
measure zero.
If $\Omega$ is a saturated sphere packing, then all Voronoi domains
are compact sets with diameter bounded by $4 \sqrt{2}$.

\begin{defi}\label{de32}
{\em The} Delaunay triangulation {\em  associated to a set $\Omega$ is
dual to the Voronoi tesselation.
It contains an edge between every pair of vertices that have Voronoi domains
that share a common face. Suppose now that the points of 
$\Omega$ are in general
position, which means that each corner of a Voronoi domain has
exactly four incident Voronoi domains.
In this case these Voronoi domains have between them four faces that
touch this corner, and these faces in turn determine (four edges of)
a Delaunay simplex. The resulting Delaunay simplices partition $\RR^3$
and make up the Delaunay triangulation. In  the
case of non-general position $\Omega$ the Delaunay triangulation is
not unique. The
possible Delaunay 
triangulations are
determined locally as limiting cases of general position points. (There are
only finitely many triangulations possible in any bounded region of space.)}
\end{defi}

%In two dimensions the optimal sphere packing density 
%$\frac {\pi}{\sqrt{12} = 0.909$ has been establshed by local
%density inequalities using either the Voronoi tesselation or
%the Delaunay triangulation.

All the simplices in a Delaunay triangulation have vertices
$\bv_i \in \Omega$ and contain no other point $\bv \in \Omega$.
We define, more generally:
\begin{defi}\label{de21}
{\rm A {\em $D$-simplex} (or 
{\em weak Delaunay-simplex})
for $\Omega$ is any tetrahedron $T$ with vertices 
$\bv_1, \bv_2, \bv_3, \bv_4 \in \Omega$ such that no other 
$\bv \in \Omega$ is in the closure of $T$.
We denote it $D( \bv_1, \bv_2, \bv_3, \bv_4)$.
}
\end{defi}

In the literature a {\em Delaunay
simplex} associated to a point set $\Omega$ is any simplex with vertices
in $\Omega$ whose circumscribing sphere
 contains no other vertex of $\Omega$ in its interior.
All simplices in a Delaunay triangulation of $\Omega$ are
Delaunay simplices, so are 
necessarily $D$-simplices, but the converse need not hold.

The admissible partitions that have been seriously studied all consist of a
domain $V(\bv )$ associated to each vertex $\bv \in V$, which we call a
{\em $V$-cell}, together with a collection of certain D- simplices 
$D(\bv_1, \bv_2, \bv_3,\bv_4)$ which we call the {\em $D$-system} of the
partition.
We use the term {\em $D$-set} to refer to a  $D$-simplex included
in the $D$-system.
We note that a $V$-cell may consist of several polyhedral pieces, 
and may even be disconnected.

The  original approach of L. Fejes-T\'oth to getting local upper 
bounds for sphere packing in $\RR^3$
used the Voronoi tesselation associated to $\Omega$.
If $\Omega$ is a saturated packing, then each Voronoi domain 
$V_{vor}( \bv )$ is a 
bounded polyhedron
consisting of points within distance at most $4 \sqrt{2}$ of $\bv$.
Examples are known of Voronoi domains in a saturated packing
that have 44 faces;
an upper bound for the number of faces of a Voronoi domain of a 
saturated packing is 49.
The {\em Voronoi partition} takes the $V$-sets
$V(\bv )$ to be the Voronoi domains of $\bv$,
with no $D$-sets,
and the vertex $D$-star $\sD_{vor} ( \bv )$ is just $V_{vor}( \bv )$.
A {\em Voronoi scoring rule} is:
\beql{302a}
Score ( \sD_{vor} (\bv )) : =( \rho (V( \bv )) -  B)~~vol(V(\bv)), ~.
\eeq
with the score constants $A = 1$ and $B$ is to be chosen optimally.
Such scoring rules are  admissible.
However it has long been known that no Voronoi scoring rule gives an optimal
inequality \eqn{M213}.
The dodecahedral conjecture states that the maximum packing density of a 
Voronoi domain is attained for a local configuration of 12 spheres 
touching at the center of faces of a circumscribed regular dodecahedron.

\paragraph{Dodecahedral Conjecture.}
{\em For a Voronoi domain $V_{vor}( \bv )$ of a unit sphere packing,
\beql{301a}
\rho (V_{vor}( \bv )) \le
\frac{\pi}{15 (1- \cos \frac{\pi}{5} \tan \frac{\pi}{3} )} \simeq
0.754697
\eeq
and equality is attained for the dodecahedral configuration.
}

A proof of the dodecahedral conjecture has been announced by
Hales and McLaughlin \cite{HM},
based on similar ideas to the Hales' approach to the Kepler conjecture.

In 1953 L. Fejes-T\'oth \cite{FT}  proposed that an optimal inequality
might exist based on a weighted averaging over Voronoi domains near a given
sphere center, and in 1964 he made a specific proposal for
such an optimal inequality.
In the notation of this paper he used the Voronoi partition and an 
admissible scoring function of the form
\beql{N303a}
\sigma (V_{vor}( \bw ), \bv ) := \omega (\bw , \bv )\{ ( A \rho (V_{vor}( \bw )) 
- B) vol ~(V_{vor}( \bw ))\},
\eeq
in which $A = 1$ and $B = \frac {\pi}{\sqrt{18}}$ and 
the {\em weights} $\omega (\bw, \bv )$ are given by
\beql{N303b}
\omega ( \bw, \bv ) := \left\{
\begin{array}{cll}
\frac{1}{12} & {\rm if} &  2  \leq  \| \bw - \bv \| \le 2 + t~, \\ [+.2in]
0 & {\rm if} & \| \bw - \bv \| > 2 + t ~,
\end{array}
\right.
\eeq
and
\beql{N303c}
\omega (\bv, \bv ) := 1 - \sum_{\bw \ne \bv} \omega (\bw, \bv).
\eeq
Here $t \ge 0$ is a fixed constant. The resulting inequality
in Theorem~\ref{th21} is optimal if the associated constant 
$\theta_{\sP,\omega}(A,B) = 0$. 
The Fejes-T\'oth scoring function corresponds to a weighted averaging over
the spheres touching a central sphere at $\bv$, where the weight assigned the
central sphere depends on how many spheres touch it.
L. Fejes-T\'oth considered choosing $t$ as large as possible 
consistent with requiring that 
$\omega (\bv, \bv ) \ge 0,$
which is equivalent to requiring that it is impossible to pack $13$ spheres
around a given sphere, with all 13 sphere centers within distance $2 + t$
of the center of the given sphere.
In 1964 he suggested \cite[p.299]{FT1} that one could take the value 
$t = 0.0534,$
 and we consider this to be Fejes-T\'oth's candidate for an
optimal inequality. It has not been demonstrated that
 $\omega (\bv, \bv ) \ge 0$ holds for this value of $t$, but note that
the argument of Theorem~\ref{th21} is valid
for any value of $t$, even if some $\omega (\bv, \bv )$ are negative. The
only issue is whether the resulting sphere packing bound is
optimal. Fejes-T\'oth \cite{FT1} explicitly noted that establishing 
an optimal inequality, if it is true,
 reduces the problem in principle to one in a finite number of variables,
possibly amenable to solution by computer.

In 1993 Wu-Yi Hsiang \cite{Hs} studied a variant of the Fejes-T\'oth
approach.
He used the Voronoi partition and an admissible scoring function
\footnote{We have converted the locally averaged density in Hsiang 
\cite[Section 3]{Hs} to the form given in \S2 by clearing denominators,
and we have cancelled out Hsiang's factor of 13. Note also that each
Voronoi domain contains exactly one sphere, so that
$\rho(V(\bv)) vol(V(\bv))= \frac {4 \pi}{3}.$}
of the form \eqn{N303a}, with the same $A=1$ and $B= \frac{\pi}{\sqrt{18}}$ 
and with the {\em weights} $\omega (\bw, \bv ) \ge 0$ given by
\beql{N304}
\omega ( \bw, \bv ) := \left\{
\begin{array}{cll}
\frac{1}{1+ N(\bv )} & {\rm if} & \| \bw - \bv \| \le \frac{218}{100} ~, \\ [+.2in]
0 & {\rm if} & \| \bw - \bv \| > \frac{218}{100} ~,
\end{array}
\right.
\eeq
where
\beql{N305}
N( \bv ) : = \# \left\{
\bw \in \Omega : ~ 0 < \| \bw - \bv \| \le
\frac{218}{100} \right\} 
\eeq
counts the number of ``near neighbors'' of $\bv$.
Hsiang announced that his local inequality is optimal
(with $\theta_{\sP,\omega}(A,B)= 0$), and that he had proved it,
which would then constitute a proof of Kepler's conjecture.
However his proof of optimality is regarded as incomplete by the mathematical
community, see G. Fejes-T\'oth's review of Hsiang's paper in 
Mathematical Reviews, the critique in \cite{Ha94}, and 
Hsiang's rejoinder~\cite{Hs95}.

In 1992 Hales (\cite{H1}, \cite{H2}) studied the 
Delaunay triangulation, which partitions $\RR^3$
into $D$-sets.
There are a finite number of 
(local) choices for 
Delaunay triangulations
of a neighborhood of a fixed $\bv \in \Omega$. 
Hales used the following function in defining his associated
weight function.

\begin{defi}\label{de25b}
{\rm
The {\em compression} of a finite region $R$ in $\RR^3$ with respect 
to a sphere packing
$\Omega$ is
\beql{201}
\Gamma (R) := (\rho ( R) - \delta_{oct}) vol(R)
\eeq
in which
\beql{202}
\delta_{oct} := \frac{-3 \pi + 12 \arccos ( \frac{1}{\sqrt{3}} )}{\sqrt{8}}
\approx 0.720903
\eeq
is the packing density of the regular octahedron of sidelength 2 
with unit spheres centered at its vertices.
%[part I, p. 3]
}
\end{defi}

Hales initially considered the admissible weight function
\beql{302}
\sigma (D(\bv_1, \bv_2, \bv_3, \bv_4 ) , \bv_i ) :=
\Gamma(D(\bv_1, \bv_2, \bv_3, \bv_4 )).
\eeq
The vertex $D$-star of $\bv$ consists of all the simplices in the
Delaunay triangulation that have $\bv$ as a vertex; we call this
set of simplices the {\em Delaunay D-star}
$\sD_{Del} (\bv )$ at $\bv$. Hales used score
constants $A=4$ and $B= 4\delta_{oct}$, with $A = 4$ used since 
each simplex is counted four times.
However he discovered that the {\em pentagonal prism} 
attained a score value exceeding what is needed to
prove Kepler's conjecture.
The pentagonal prism is conjectured to be extremal for this score function.

The fact that the (conjectured) extremal configurations for the 
Voronoi tesselation and Delaunay triangulation do not coincide suggested 
to Hales that a 
hybrid scoring rule be considered that combines the best features 
of the Voronoi and Delaunay scoring function.
In 1997 Hales again considered a 
Delaunay triangulation, but modified
the scoring rule to depend on the shape of the D-simplex 
$D( \bv_1, \bv_2, \bv_3, \bv_4 )$. For some simplices he used the
weight function above, while for others he cut the simplex into
four pieces, one for each vertex, call the pieces $V(D, \bv_i),$
and assigned the weights\footnote{More precisely, he used the
``analytic continuation'' of this scoring function that is described
in Appendix A.}
$$\sigma(D(\bv_1, \bv_2, \bv_3, \bv_4 ), \bv_i)= 4\Gamma( V(D, \bv_i)), $$
for $1 \leq i \leq 4.$ 
He also partitioned
a vertex D-star into pieces called
``clusters'' whose score functions could be evaluated separately
and added up to get the total score. Each ``cluster'' is a
finite union of Delaunay simplices filling up that part of the
vertex D-star at $\bv$ lying in a pointed cone with vertex $\bv$.
This vertex cone subdivision
facilitates computer-aided proofs by decomposing the problem
into smaller subproblems.
Hales (\cite{I}, \cite{II}) presented evidence that this 
modified scoring function 
satisfies an optimal local inequality.
He showed that the two known local extremal configurations 
\footnote{These correspond to Voronoi cells being the rhombic
dodecahedron or trapezo-rhombic dodecahedron in Fejes-T\'oth 
\cite[p. 295]{FT1}.} gave 
local maxima of the score of the Delaunay $D$-star in the 
configuration space with
\beql{303}
Score(\sD_{Del} (\bv )) = 8 pt~,
\eeq
where
\beql{303a}
pt := \frac{11 \pi}{3} - 12 \arccos \left(
\frac{1}{\sqrt{3}} \right) \simeq 0.0553736 ~,
\eeq
which is the optimal value.
He also showed that this was a global upper bound over the subset of 
 configurations  described by a vertex map\footnote{See \S5 for a 
definition of vertex map $\sG (\bv)$.}
$\sG ( \bv )$ that is triangulated.
Hales \cite[Conjecture 2.2]{I} conjectured that the 
modified score function achieved the optimal inequality
\beql{304}
Score (\sD_{Del} ( \bv )) \le 8 pt ~,
\eeq
for all Delaunay D-stars $\sD_{Del} (\bv )$.
However he and his student S. P. Ferguson \cite{FH} discovered  
that a pentagonal prism configuration comes very close to 
violating the inequality \eqn{304}.
Furthermore there turned out to be many similar difficult configurations 
which might possibly violate the inequality.
These and other difficulties indicated that it was not 
numerically feasible to prove \eqn{304}
by a computer proof, assuming  that \eqn{304} is actually true.

Hales and Ferguson together then
further modified both the partition rule $\sP$ and 
the scoring rule $\sigma$, to obtain a rule
with the following properties.
\begin{itemize}
\item[(i)]
It makes the  score inequality stronger on the known bad cases related to the
pentagonal prism configuration.
\item[(ii)]
It uses a more complicated notion of ``cluster'', which includes
Voronoi pieces as well as D-sets, and which retains 
the ``decoupling'' property that it
is completely determined by 
vertices of $\Omega$ in the cone above it.
\item[(iii)]
It chooses a scoring function which when combined with
``truncation'  on clusters is still strong enough to rule out 
most configurations. The ``truncation'' operation greatly reduces 
the number of configurations to be checked, at the cost of weakening
the inequality to be proved.
\end{itemize}
In \S4 we give a precise description of the Hales-Ferguson rules 
$(\sP , \sigma )$.
\section{Hales-Ferguson Partition and Score Function}
\hsp
Ferguson and Hales \cite{FH} use the following partition 
and scoring rule.
The partition uses two types of $D$-simplices, with a complicated 
rule for picking which
ones to include as $D$-sets in the partition.
Modified Voronoi domains $V(\bv )$ are used as $V$-sets.
These differ from the usual Voronoi domain (with the $D$-sets removed) 
by mutually exchanging some regions called ``tips''.
The scoring rule is also complicated:
the weight function used on a $D$-simplex no longer depends on 
just its shape,
but depends on the structure of nearby $D$-sets.

We begin by defining the two types of $D$-simplices.
\begin{defi}\label{de41}
A QR-tetrahedron {\rm (or} quasi-regular tetrahedron$)$
{\rm is any tetrahedron with all vertices in $\Omega$ and all edges of 
length $\le \frac{251}{100}$.
}
\end{defi}

\begin{defi}\label{de42}
{\rm A} QL-tetrahedron {\rm (or} quarter$)$
{\rm is any tetrahedron with all vertices in $\Omega$ and five edges of 
length
$\le \frac{251}{100}$ and one edge with length 
$\frac{251}{100} < l \le 2 \sqrt{2}$.
The long edge is called the {\em spine} (or {\em diagonal}) of the
QL-tetrahedron.
%It is {\em strict} if $l > \frac{251}{100}$.
}
\end{defi}

For some purposes\footnote{In proving inequalities one wants to work on
a compact set. In compactifying the space of
configurations, this requires allowing the lower inequality in the
definition of $QL$-tetrahedron to be an equality.} 
the case of a spine of length exactly $\frac{251}{100}$ should be considered 
as either a QL-tetrahedron and QR-tetrahedron.
Here we treat it exclusively as a QR-tetrahedron.

Neither kind of tetrahedron is guaranteed to be included in the Delaunay 
triangulation of $\Omega$,
but we do have:
\begin{lemma}\label{le41}
All QR-tetrahedra and QL-tetrahedra are D-simplices.
\end{lemma}

The Hales-Ferguson partition rule starts by selecting which
$D$-simplices to include in the $D$-system.
These consist of:
\begin{itemize}
\item[(i)]
All QR-tetrahedra.
\item[(ii)]
Some QL-tetrahedra. The QL-tetrahedra included in the partition
 satisfy the {\em common spine condition} which states 
that for a given spine, either all QL-tetrahedra having that spine are 
included, or none are.
\end{itemize}
This collection of tetrahedra must (by definition of admissible partition) 
form a nonoverlapping set,
where we say that two sets $S_1$ and $S_2$ {\em overlap} if 
$\bar{S_1} \cap \bar{S_2}$ has positive Lebesgue measure in $\RR^3$.
To justify (i) we have:
\begin{lemma}\label{le42}
No two QR-tetrahedra overlap.
\end{lemma}

$QL$-tetrahedra may overlap $QR$-tetrahedra or other $QL$-tetrahedra, 
hence one needs a rule
for deciding which $QL$-tetrahedra to include.
To begin with, $QL$-tetrahedra can overlap $QR$-tetrahedra in 
essentially one way.

\begin{lemma}\label{le43}
If a $QL$-tetrahedron and $QR$-tetrahedron overlap, then the $QR$-tetrahedron
has a common face with an adjacent $QR$-tetrahedron, and the two unshared
vertices of these $QR$-tetrahedra are the endpoints of the spine of the
$QL$-tetrahedron. 
The union of these two $QR$-tetrahedra can be
partitioned into three $QL$-tetrahedra having the given spine, 
which includes the given $QL$-tetrahedron. Aside from these
$QL$-tetrahedra, no other $QL$-tetrahedron
overlaps either of these two $QR$-tetrahedra.
\end{lemma}

This lemma shows that the $QL$-tetrahedra having a given spine
have the property that either all of them or none of them overlap 
the set of  $QR$-tetrahedra. We next consider how $QL$-tetrahedra can overlap
other $QL$-tetrahedra.  The following configuration plays an
important role.

\begin{defi}\label{Nde43}
{\rm A {\em $Q$-octahedron} is an octahedron whose 6 vertices
$\bv_i \in \Omega$ and whose 12 edges each have lengths 
$2 \le l \le \frac{251}{100}$.}
\end{defi}
 
A $Q$-octahedron has three interior diagonals.
If a diagonal has length $2 < l \le \frac{251}{100}$ then it partitions the 
$Q$-octahedron into four $QR$-tetrahedra.
If a diagonal has length $\frac{251}{100} < l \le 2 \sqrt{2}$ then it partitions 
the $Q$-octahedron into four $QL$-tetrahedra of which is the common spine.
If a diagonal has length $l > 2\sqrt{2}$ it yields no partition.
A $Q$-octahedron thus gives between zero and three
different partitions into 
four $QR$-tetrahedra or $QL$-octahedra.
We call it a {\em live $Q$-octahedron} if it has at least one such partition.
Lemma~\ref{le43} implies that if it has a partition into $QR$-octahedra
then it has no other partition into $QR$-octahedra or $QL$-octahedra.

\begin{lemma}\label{Nle43}
A $QL$-tetrahedron having spine a diagonal of a $Q$-octahedron 
does not overlap any $QL$-tetrahedron whose spine is not a diagonal 
of the same $Q$-octahedron.
\end{lemma}

The rules for choosing which $QL$-tetrahedra to include either take all 
$QL$-tetrahedra having the same spine or take none of them.
Thus the selection rule really specifies which spines to include.
\begin{defi}\label{Nde44}
{\rm Consider an edge $[\bv_1, \bv_2]$ with $\bv_1, \bv_2 \in \Omega$ and
$\frac{251}{100} < \| \bv_1 - \bv_2 \| \le 2 \sqrt{2}$.
A vertex $\bw \in \Omega$ is called an {\em anchor} of the edge
$[\bv_1 , \bv_2 ]$ if
$$\| \bw- \bv_i \| \le \frac{251}{100} \quad\mbox{for}\quad
i=1,2 ~.
$$
}
\end{defi}

Ferguson and Hales  use the number of $QL$-tetrahedra having a 
spine $[\bv_1, \bv_2 ]$ 
and the number of anchors of that spine in deciding which $QL$-tetrahedra 
to include in the $D$-system.
Call a $QL$-tetrahedron {\em isolated} if it is the only $QL$-tetrahedron 
on its
spine $[\bv_1, \bv_2 ]$.
The inclusion rule for an isolated
$QL$-tetrahedron is:
\begin{itemize}
\item[(QL0)]
An isolated $QL$-tetrahedron is included in the $D$-system if and only 
if it
overlaps\footnote{It cannot overlap a $QR$-tetrahedron by 
Lemma \ref{le43}.} no other $QL$-tetrahedron or $QR$-tetrahedron.
\end{itemize}
Next consider spines $[\bv_1, \bv_2 ]$ which have two or more 
associated $QL$-tetrahedra.
Such spines have at least three anchors, and the inclusion rules are:
\begin{itemize}
\item[(QL1)]
Each non-isolated $QL$-tetrahedron on a spine with 5 or more anchors
is included in the $D$-system.
\item[(QL2)]
Each non-isolated $QL$-tetrahedron on a spine with 4 anchors is 
included in the $D$-system,
if the spine is not a diagonal of some $Q$-octahedron.
In the case of a live $Q$-octahedron, we include all $QL$-tetrahedra having one
particular diagonal. and exclude 
all $QL$-tetrahedra on other diagonals.
For definiteness, we choose the spine to 
be the shortest diagonal. In case of a tie for shortest diagonal, a
suitable tie-breaking rule is used. 
\item[(QL3)]
Each non-isolated $QL$-tetrahedron on a spine with 3 anchors is 
included in the $D$-system if each $QL$-tetrahedron on the spine
does not overlap any other $QL$-tetrahedron or $QR$-tetrahedron,
or overlaps only isolated $QL$-tetrahedra.
It is excluded from the $D$-system if some tetrahedron on the spine
overlaps either a $QR$-tetrahedron
or a non-isolated $QL$-tetrahedron having four or more anchors.
Finally, if some tetrahedron on the spine overlaps a 
nonisolated $QL$-tetrahedron having exactly
three anchors, then the spine of the overlapped set is unique,
and exactly one of these two sets of non-isolated
$QL$-tetrahedra with 3 anchors is to be included in the $D$-system,
according to a tie-breaking rule\footnote{The tie-breaking rule could
be to include the spine with lowest endpoint using a lexicographic
ordering of points in $\RR^3$. It appears to me that
Hales would permit an
arbitrary choice of which one to include, see Lemma~\ref{Nle44} (iii).}.
\end{itemize}
The set of $QL$-tetrahedra selected above are pairwise disjoint and 
are disjoint from all $QR$-tetrahedra.
This is justified by the following lemma.
\begin{lemma}
\label{Nle44}
(i) If two $QL$-tetrahedra overlap, then at most one of them has $5$ or 
more anchors.

(ii)~If two overlapping $QL$-tetrahedra each have $4$ anchors, then their 
spines are (distinct) diagonals of some $Q$-octahedron.

(iii)~If a nonisolated $QL$-tetrahedron with $3$
anchors overlaps another $QL$-tetrahedron having three anchors, 
then each of their spines contains exactly two nonisolated $QL$-tetrahedra,
and these four $QL$-tetrahedra overlap no other $QL$-tetrahedron or
$QR$-tetrahedron.
\end{lemma}

We call the set of $QR$-tetrahedra and $QL$-tetrahedra selected as above the 
{\em {\em Hales-Ferguson} $D$-system.} 
(Hales and Ferguson call this a  $Q$-system.)

We now define the $V$-cells of the Hales-Ferguson partition.
To begin with, we take the Voronoi domain $V_{vor} (\bv )$ at vertex $\bv$
and remove from it all $D$-simplices in the $D$-system to obtain a reduced
Voronoi region $V_{red} ( \bv )$.
Next we move certain regions of $V_{red} (\bv )$ called ``tips'' to
neighboring reduced Voronoi regions to obtain modified regions
$V_{mod} (\bv )$ and finally we define the $V$-cell $V( \bv )$ at $\bv$ to be
the closure of $V_{mod} (\bv )$.

\begin{defi}\label{Nde45}
{\rm Let $T$
be any tetrahedron such that the center $ \bx = \bx (T)$ of its
circumscribing sphere lies outside $T$.
A vertex $\bv$ of $T$ is {\em negative} if the plane $H$ determined
by the face $F$ of $T$ opposite $\bv$ separates $\bv$ from $\bx$.
The ``{\em tip}'' $\Delta (T,\bv )$ of $T$ associated to a negative 
vertex $\bv$
is that part of the Voronoi region of $\bv$ with respect to the points
$\{ \bv_1, \bv_2, \bv_3, \bv_4 \}$ that lies in the closed half
plane $H^+$ determined by $H$ that contains $\bx$.
The ``tip'' region $\Delta (T,\bv )$ does not overlap $T$, and
is a tetrahedron having $\bx$ as a vertex,
and has three other vertices lying on $H$.
}
\end{defi}
\begin{lemma}\label{Nle45}
(i) A $QR$-tetrahedron or $QL$-tetrahedron $T$ has at most
one negative vertex.

(ii) If a negative vertex is present, then the three vertices of
the associated ``tip'' that lie on $H$ actually lie in 
the face of $T$ opposite to the negative vertex.

(iii)~The ``tip'' of any tetrahedron in the Hales-Ferguson $D$-system
either does not overlap any $D$-simplex in the Hales-Ferguson $D$-system,
or else is entirely contained in the union of the $D$-simplices in the 
$D$-system.
\end{lemma}

We say that a ``tip'' that does not overlap any $D$-set is {\em uncovered.}
The lemma shows that uncovered ``tips'' lie in the union of the
Voronoi regions $\{V_{red} (\bw ) : \bw \in \Omega \}$, so that 
rearrangement of uncovered ``tips'' is legal.
There is an a priori possibility that two ``tips'' may overlap\footnote{
I don't know if this possibility can occur.} each other.

\paragraph{Uncovered Tip Rearrangement Rule.}
Each $\by \in \RR^3$ that belongs to an uncovered  ``tip'' is reassigned to 
the nearest {\em vertex} $\bw \in \Omega$ such that
$\by$ is not in an uncovered ``tip'' of 
 any pair $(T, \bw )$ where $T$ is in the Hales-Ferguson $D$-system 
and $\bw$ is a negative vertex of $\bv$.
(A tiebreaking rule is used if two nearest vertices $\bw$ are equidistant.)

This rule cuts an uncovered ``tip'' into a finite number of polyhedral pieces 
and reassigns the pieces to different reduced Voronoi regions.
This prescribes how $V_{mod} (\bv )$ is constructed, and thus defines the
Hales-Ferguson $V$-cells $V(\bv )$.

\begin{figure}[htb]
\begin{center}
\input tip2.pstex_t
\end{center}
\caption{``Tip'' of $D$-simplex $[\bv_1 , \bv_2, \bv_3 ]$}
\label{fg101}
\end{figure}

A two-dimensional analogue of a ``tip'' is pictured\footnote{
See Figure 2.1 of Hales \cite{II} for another example.} in Figure 1.
In this figure the  triangle $T = [ \bv_1, \bv_2, \bv_3 ]$ plays the role of a 
$D$-simplex, with $\bv_2$ as a  negative vertex and the ``tip'' is the
shaded region. The points $\bc_{012}, \bc_{013}, \bc_{123}$ are centroids
of the triangles determined by the  corresponding $\bv_i$'s.
The shaded triangle $[ \bc_{012}, \bc_{013}, \bc_{123}]$ is in the
Voronoi cell $V_{vor}(\bv_0)$ while the remainder of the ``tip'' is
in the Voronoi cell $V_{vor}(\bv_2).$ The uncovered tip rearrangement rule
partitions the part in $V_{vor}(\bv_2)$ into three triangles which
are reassigned to the $V$-cells $V(\bv_0), V(\bv_1)$ and $ V(\bv_3),$
e.g. $[ \by_1, \by_2, \bc_{012}]$ is reassigned to  $ V(\bv_1)$.
The reassignment of the ``tip'' 
ensures that the pointed cone over $\bv_2$
generated  by the $D$-simplex $[\bv_1, \bv_2, \bv_3]$ does not
contain any part of the $V$-cell at $\bv_2$. In this example 
the $V$-cell at $\bv_0$ does not feel the effect 
of the vertex $\bv_2$, due to the rearrangement.

We now turn to the Hales-Ferguson scoring rules.
These use the compression function $\Gamma (S)$ given in \eqn{201}.
The compression function is additive:
If $S= S_1 \cup S_2$ is a partition, then
\beql{203}
\Gamma (S) = \Gamma (S_1 ) + \Gamma (S_2 ) ~.
\eeq
For a $D$-simplex $T$,
\beql{204}
vol (T) \rho (T) = \sum_{i=1}^4
\frac{(\mbox{solid angle})_i}{3} ~,
\eeq
where a full solid angle is $4 \pi$.

The Hales-Ferguson weight function for a V-cell is as follows.
\begin{itemize}
\item[(S1)]
For a $V$-cell $V(\bv )$,
\beql{N41}
\sigma_{HF} (V( \bv ), \bw ) = \left\{
\begin{array}{cll}
4 \Gamma ( V( \bv )) & \mbox{if} & \bv = \bw ~, \\ [+.2in]
0 & \mbox{if} & \bv \neq \bw ~.
\end{array}
\right.
\eeq
\end{itemize}

We next consider the weight function for $D$-sets. 
Let $(T, \bv )$ denote a $D$-simplex together with
a vertex $\bv$ of it. 
\begin{defi}\label{Nde46}
{\rm The {\em Voronoi measure} $vor (T, \bv )$ is defined as follows.
If the center of the circumscribing sphere of $T$ lies inside
$T$, then $T$ is partitioned into four pieces
$$V_{vor}^+ (T, \bv_i ) := \{
\bx \in T: \| \bx - \bv_i \| \le \| \bx - \bv_j \| \quad\mbox{for}\quad
1 \le j \le 4 \} 
$$
and then
\beql{N42}
vor (T, \bv ) := \Gamma ( V_{vor}^+ (T, \bv )) ~.
\eeq
There is an analytic formula for the right side of \eqn{N42} given in 
Appendix A, and this formula is used to define $vor (T, \bv )$ in  
cases where the circumcenter falls outside $T$.
}\end{defi}

In cases where the circumcenter is outside $\bv$, and $\bv$ is a 
negative vertex, then
\beql{N43}
vor (T, \bv ) = \Gamma ( V_{vor} (T, \bv ) \cup \mbox{``tip''} )
\eeq
while for the other three vertices parts of the ``tip'' 
are counted with a negative
weight, in such a way that
\beql{N44}
\sum_{i=1}^4 vor (T, \bv_i ) = 4 \Gamma ( T)
\eeq
holds in all cases. The weight function for a $D$-set is given as
follows:
\begin{itemize}
\item[(S2)]
For a $QR$-tetrahedron $T= D(\bv_1, \bv_2, \bv_3, \bv_4 )$ in 
the $D$-system,
\beql{N45}
\sigma_{HF} (T, \bv ) =
\left\{
\begin{array}{cl}
\Gamma (T) & \mbox{if the circumradius of $T$ is at most 
$\frac{141}{100}$.} \\ [+.2in]
vor (T, \bv ) & \mbox{if the circumradius of $T$ exceeds $\frac{141}{100}$.}
\end{array}
\right.
\eeq
\end{itemize}
The $QL$-tetrahedron scoring function is complicated. For a
$QL$-tetrahedron $T$. let $\eta^+ (T)$ be the maximum of the circumradii 
of the two 
triangular faces of $T$ adjacent
to the spine of $T$, and define the function
\beql{N46}
\mu (T, \bv ) :=
\left\{ \begin{array}{cll}
\Gamma (T) & \mbox{if} & \eta^+ (T) \le \sqrt{2} \\ [+.2in]
vor (T, \bv ) & \mbox{if} & \eta^+ (T) > \sqrt{2}
\end{array}
\right.
\eeq
Then the $QL$-tetrahedron scoring function is defined by:
\begin{itemize}
\item[(S3)]
(``Flat quarter'' case)
For a $QL$-tetrahedron $T$ and a vertex $\bv$ not on its spine,
\beql{N47}
\sigma_{HF} (T, \bv ) := \mu (T, \bv ) ~.
\eeq
\item[(S4)]
(``Upright quarter case'')
For a $QL$-tetrahedron $T$ with vertex $\bv$ on its spine, let
$\hat{\bv}$ denote the opposite vertex on the spine.
If $T$ is an isolated $QL$-tetrahedron, set
\beql{N48}
\sigma_{HF} (T, \bv ) := \mu (T, \bv ) ~.
\eeq
If $T$ is part of a $Q$-octahedron, set
\beql{N49}
\sigma_{HF} (T, \bv ) := \frac{1}{2} (
\mu (T, \bv ) + \mu (T, \hat{\bv} )) ~.
\eeq
In all other cases, set
\beql{N410}
\sigma_{HF} (T, \bv ) := \frac{1}{2}
( \mu (T, \bv ) + \mu (T, \hat{\bv} )) +
\frac{1}{2} ( vor_0 (T, \bv ) - vor_0 (T, \hat{\bv} )) ~,
\eeq
in which $vor_0 (T, \bv )$ is a ``truncated Voronoi measure''
that only counts volume within radius $\frac{1}{2} ( \frac{251}{100} )$ of 
vertex $\bv$, which is defined in Appendix A, and in \cite[pp. 9-11]{FH}.
\end{itemize}
The scoring rule (S4) is the most complicated one. In it
the definition \eqn{N49} plays an important role in 
obtaining good  bounds
for the pentagonal prism case treated in Ferguson \cite{Fe}, while
the definition \eqn{N410} is important in analyzing general configurations
using truncation in Hales \cite{IV}.

\begin{theorem}\label{Nth41}
The Hales-Ferguson partition and scoring function
$(\sP_{HF} , \sigma_{HF} )$ are admissible, with score constants $A = 4$ and
$B = 4 \delta_{oct}.$
\end{theorem}

\paragraph{Proof.}
It is easy to verify that the 
definitions for scoring
$QR$-tetrahedra and $QL$-tetrahedra satisfy the weighted density average
property
\beql{N411}
\sum_{i=1}^4 \sigma_{HF} (T, \bv_i ) = 4 \Gamma (T) ~,
\eeq
which correspond to $A=4$ and $B = 4 \delta_{oct},$ using \eqn{N44}.
Most of the remaining admissibility conditions are verified by 
Lemmas \ref{le41}--\ref{Nle45} except for locality.
For locality, a  conservative estimate indicates that the rules for 
removing and adding
``tips'' to determine the $V$-cell $V( \bv )$ are determined by sphere centers
$\bw \in \Omega$ with $\| \bw - \bv \| \le 12 \sqrt{2}$.
Finally the  score function on the $D$-simplices is determined by vertices
within distance $6 \sqrt{2}$ of $\bv$.~~~$\bsq$

Theorem \ref{th21} associates to $(\sP_{HF}, \sigma_{HF} )$ a sphere-packing
bound that the Hales program asserts is optimal.
To establish the Kepler bound
\beql{213}
\bar{\rho} (\Omega) \le \frac{\pi}{\sqrt{18}} ~,
\eeq
via \eqn{M209}, one must prove that
\beql{214}
 \theta :=\theta_{\sP_{HF}, \sigma_{HF}}(4, 4 \delta_{oct}) = 8 pt~,
\eeq
where
\beql{215}
pt :=
\frac{11 \pi}{3} - 12 \arccos \left( \frac{1}{\sqrt{3}} \right) \simeq
0.0553736 ~.
\eeq

The score function
$Score(\sD_{HF}({\bv}))$ is discontinuous as a function of the sphere centers
in $\Omega$ near $\bv$, because it
is a sum of contributions of pieces which may appear and disappear
as sphere centers move, and discontinuities occur when $QL$-tetrahedra 
convert to $QR$-tetrahedra.
To deal with this, one compactifies the configuration space 
 by allowing some sphere center configurations
to have more than one legal decomposition into pieces (but at most 
finitely many).
The optimization problem can then be split into a finite number of subproblems
on each of which $\sigma_{HF}$ is continuous.

The complexity of the definition of $(\sP_{HF}, \sigma_{HF} )$ is designed
to yield a computationally tractable nonlinear
optimization problem.
The introduction of $QL$-tetrahedra and the complicated score function on 
them is designed to help get good bounds for the pentagonal prism case and
similar cases.
The rule for moving ``tips'' is intended to facilitate decomposition of the 
nonlinear optimization problem into more tractable pieces via 
Theorem \ref{le55} below, and the use of  ``truncation.''

\section{Kepler Conjecture}
\hsp
The main result to be established by the Hales program is the following.

\begin{theorem}\label{th51}
{\rm (Main Theorem)}
For the Hales-Ferguson partition and scoring rule 
$(\sP_{HF} , \sigma_{HF} )$, and any
$\bv \in \Omega$ in a saturated sphere-packing, the vertex $D$-star
$\sD_{HF} (\bv )$ at $\bv$ satisfies
\beql{501}
Score (\sD_{HF} ( \bv )) \le 8 pt ~,
\eeq
where 
$pt := \frac{11 \pi}{3} - 12 \arccos \left( \frac{1}{\sqrt{3}} \right) \simeq 0.0553736$.
\end{theorem}

The Kepler conjecture follows by Theorem \ref{th21}.

To prove the inequality \eqn{501}, by translation-invariance we can reduce 
to the case $\bv = \bzero$ and search the set of all possible vertex stars, 
which by
\S4 are determined by those points $\bw \in \Omega$ with 
$\| \bw \| \le 12 \sqrt{2}$. From now on we assume 
$\bzero \in \Omega$ and $\bv = \bzero$.

The space of possible sphere centers 
$\{ \bw \in \Omega: \| \bw \| \le 12 \sqrt{2} \}$ is compact.
It can be decomposed into a large number of pieces, on 
each of which the score function is continuous.
To obtain {\em compact} pieces, we must compactify the
configuration space  by assigning more than one possible 
local D-star $\sD (\bv )$
to certain arrangements  
of sphere centers. The compactification assigns at most finitely 
many to each arrangement
with an absolute upper bound on the number of possibilities.

The definition of the score $Score(\sD ( \bv ))$ involves a sum over
the $V$-sets and $D$-sets.
The usefulness of the compression measure $\Gamma (S)$ is justified 
by the following lemma.

\begin{lemma}\label{le51}
(i) Every $QR$-tetrahedron $T$ satisfies
\beql{502}
\Gamma (T) \le pt ~,
\eeq
with equality occurring only when $T$ is a regular tetrahedron of 
edge length $2$.

(ii) A $QL$-tetrahedron $T$ has
\beql{503}
\Gamma (T) \le 0 ~,
\eeq
with equality occurring for those $T$ having five edges of length $2$ and 
a spine of length $2 \sqrt{2}$.
\end{lemma}

Result (ii) illustrates a somewhat counterintuitive behavior of the 
local density function:
when holding five edges of a tetrahedron fixed of length 2, and
allowing the sixth edge to vary over 
$\frac {251}{100} \leq l \leq 2 \sqrt {2}$,
 the local density measure is largest for a
spine of {\em maximal} length.

The vertices $\bw \in \Omega$ with $\| \bw \| \le \frac{251}{100}$ 
play a particularly
important role, for they determine all $QR$-simplices of $\Omega$ 
containing $\bzero$
as a vertex.

\begin{defi}\label{de51}
{\rm The {\em planar map} (or {\em graph}) $\sG ( \bv )$ associated 
to a vertex
$\bv \in \Omega$
consists of the radial projection onto the unit sphere
$\partial B (\bv; 1) = \{ \bx \in \RR^3 : \| \bx - \bv \| =1 \}$
centered at $\bv$
of all vertices $\bw \in \Omega$ with 
$\| \bw - \bv \| \le \frac{251}{100}$ plus all those edges
$[\bw, \bw']$ between two such vertices which have length
$\| \bw - \bw' \| \le \frac{251}{100}$.
}
\end{defi}

Here we regard the planar map $\sG ( \bv )$ as being given with its embedding 
as a set of arcs on the
sphere.
The following lemma asserts that no new vertices are introduced other 
than those coming from points of $\Omega$ with $\| \bw \| \le \frac{251}{100}$.

\begin{lemma}\label{le52}
The radial projection of two edges $[\bw_1, \bw_2]$, $[\bw'_1 , \bw'_2 ]$ 
as above onto the unit sphere $\partial B(\bv ;1)$ give two arcs in 
$\sG (\bv )$ which either are disjoint
or which intersect at an endpoint of both arcs.
\end{lemma}

We study local configurations classified by the planar map $\sG ( \bzero )$.
The planar map $\sG ( \bzero )$, which is determined by the vertices 
$\| \bw \| \le \frac{251}{100}$, does
not in general uniquely determine the vertex D-star $\sD_{HF} ( \bzero )$, 
but does determine all points $\bx$ in it with 
$\| \bx \| \le \frac{251}{200}$.
%In many subcases, however, the score of the truncated region yields 
%sufficiently good bounds to 
%prove the optimal inequality.

\begin{defi}
\label{deP52}
{\rm The part of $\sD_{HF} ( \bv )$ that lies in the pointed cone with 
base point $\bv$ determined by a face of the map
$\sG ( \bv )$ is called the {\em cluster} over that face. Note that
the face need not be convex, or even simply connected- it could be
topologically an annulus, for example.
}
\end{defi}

The following lemma shows that the vertex D-star
$\sD_{HF}({\bf 0})$ can be cut up into clusters in a way compatible
with the scoring function.

\begin{lemma}\label{le53}
Each $QR$-tetrahedron or $QL$-tetrahedron in the D-star $\sD_{HF}({\bf 0})$
is contained in a single cluster. Furthermore all such tetrahedra having a
common spine are contained in a single cluster.
\end{lemma}

In effect the partition of the vertex D-star into clusters 
partitions the $V$-cell into smaller pieces, while leaving the $D$-sets
unaffected. The scoring function is additive over any partition of
a $V$-cell into smaller pieces, according to \eqn{203} and \eqn{N41}.
The {\em score} $\sigma_{HF}(F)$ of the cluster determined
by a  face $F$ of $\sG ( \bv )$ is the sum of the scores of the
$QR$-tetrahedra and $QL$-tetrahedra in the cluster, plus the
Voronoi score $4\Gamma(R)$ of the remaining part $R$ of the cluster. 
We then have
\beql{503a}
Score(\sD_{HF}(\bv)) = \sum_{F \in \sG ( \bv )}\sigma_{HF}(F).
\eeq

We now consider clusters associated to the simplest faces $F$ in
 the graph $\sG ( \bv )$. Each triangular face corresonds to a 
$QR$-tetrahedron in $\sD_{HF} ( \bv )$, and, conversely, 
each $QR$-tetrahedron in $\sD_{HF} ( \bv )$ produces a triangular face. 
A {\em quad cluster} is a cluster over a quadrateral face.
A $Q$-octahedron with spine ending at ${\bf 0}$ results in a quadrilateral
face, but there are many other kinds of quad clusters. 
In the case of faces $F$ with $\ge 5$ edges, the cluster may consist of
a $V$-cell plus some $QL$-tetrahedra, in many possible ways.
All the possible decompositions into such pieces have to be 
considered as separate configurations.

\begin{lemma}\label{le54}
(i) A cluster over a triangular face $F$ consists of a single
 QR-tetrahedron, and
conversely. The score of such a cluster is at most
1 pt, and equality holds if and only if it is a regular tetrahedron
of edge length 2.

(ii) The sum of the score functions over any quad cluster is at most
zero. Equality can occur only if the four sphere centers $\bv_i$ corresponding
to the vertices of the quad cluster each
lie at distance 2 from $\bv$ and also from each other, if they
share an edge of the quad cluster.

(iii) The score of a cluster over any face with five or more
sides is strictly negative.
\end{lemma}

The extremal graphs where equality is known to occur in \eqn{501} have 
eight triangular faces and six quadrilateral faces. The upper bound
of $8 pt$ for these cases is implied by this lemma.   
(It appeared first in  Hales~\cite[Theorem 4.1]{II}.)

The following result rules out graphs $\sG ( \bv )$ with 
faces of high degree.

\begin{theorem}\label{th52}
All decomposition stars $\sD_{HF} ( \bzero )$ with planar maps 
$\sG ( \bzero )$ satisfy
\beql{504}
Score( \sD_{HF}( \bzero )) \le 8 pt
\eeq
unless the planar map $\sG (0)$ consists entirely of
(not necessarily convex) faces of the following kinds:
polygons having at most $8$ sides,
in which pentagons and hexagons may contain an 
isolated interior vertex or a single edge
from an interior vertex to an outside vertex, and a pentagon 
may exclude from its interior a
triangle with two interior vertices.
\end{theorem}

There remain a finite set of possible map structures that 
satisfy the conditions of Theorem \ref{th52}.
Here we use the fact that there can be at most 50 vertices $\bv$ with
$\| \bv \| \le \frac{251}{100}$.
The list is further pruned by various methods, and reduced to about 5000
cases. Since the (putative) extremal cases are already covered
by Lemma~\ref{le54}, in the remaining cases
one wishes to prove a strict inequality
in \eqn{501}, and such bounds can be
obtained in principle by computer.

Most of the remaining cases are eliminated by linear programming bounds.
The linear programs involve obtain upper bounds for the score function
$Score (\sD_{HF} (\bzero ))$ for a planar map $\sG$ of a particular 
configuration type, using as objective function the score function,
in the form:
\beql{P55}
Maximize ~~ Score (\sD_{HF} (0)) := \sum_{faces \atop F} \sigma (F) ~,
\eeq
where the variable $\sigma (F)$ is the sum of weights
associated to the cluster over the face $F$.
The use of linear programming relaxations of the nonlinear program
 seems to be a necessity in bounding the score function.
For example the compression function $\Gamma (R)$ for different 
regions $R$ is badly
behaved: it is neither convex nor concave in general.
The linear constraints include hyperplanes bounding the convex hull of the 
score function over the variable space.

One can decouple the contributions of the separate faces $F$ of 
$\sG (\bv )$ using the following result.
\begin{lemma}\label{le55}
{\rm (Decoupling Lemma)}
Let $\bv \in \Omega$ be a vertex of a saturated packing and let $F$ 
be a face of the associated
planar map $\sG ( \bv )$, and let $\sC_F$ denote the $($closed$)$ 
pointed cone over $F$ with vertex $\bv$, and let $\sC_{F,red}$ denote
the closure of the cone over $F$ obtained by removing from  $\sC_F$
all cones over $D$-sets with a corner at $\bv$.
Then the portion of the $V$-cell $V(\bv )$ that lies in $\sC_F$ is 
completely determined by the vertices of
$\Omega$ that fall in the smallest closed convex cone
$\bar{\sC}_F$ containing $\sC_F$.
In particular,
\beql{505}
V_F :=V( \bv ) \cap \sC_F = V(\Omega \cap \bar{\sC}_F , \bv ) \cap \sC_{F,red} ~.
\eeq
\end{lemma}

To obtain such a decoupling lemma requires the exchange  of ``tips'' 
between Voronoi domains, as described in \S4. 

The decoupling lemma permits the score function 
$\sigma (V( \bv ), \bv )$ to be decomposed into polyhedral pieces 
that depend on only a few of the nearby vertices.
This decomposes the problem into a sum of smaller problems, to
bound the scores of the pieces
$\sigma ( V( \bv ) \cap \sC_F , \bv )$ in terms of these vertices.
It will often be applied when the face $F$ is convex, 
in which case $\bar{\sC}_F = \sC_F$.

A futher very important relaxation of the linear programs involves
``truncation.'' The {\em truncated V-cell} is
\beql{505a}
V_{trunc}(\bv) := V(\bv) \cap \bB(\bv: \frac {251}{200}).
\eeq
We may consider truncation of that part of the
 $V-cell$ over each face of $\sG$ separately.

\begin{lemma}\label{le56}
Let $F$ be a face of $\sG(\bv)$ and $\sC_F$ the cone over that face.
The region $V_{trunc}(\bv) \cap \sC_F$ is entirely determined by
the vertices of $\sG(\bv)$ in $\sC_F$. If $\sV_F$ denotes this set
of vertices, together with $\bv$ then this region is the closure of
$(V_{vor}(\sV, \bv) \cap \sC_F) - \{ D-sets \}.$ The compression
function satisfies the bound
\beql{505b}
\Gamma ( V_{trunc}(\bv) \cap \sC_F) \ge \Gamma( V(\bv) \cap \sC_F).
\eeq
\end{lemma}

The inequality \eqn{505b} implies that replacing a Voronoi-type region
 by a truncated region can only increase the score, hence one can relax
the linear program
by using the score of truncated regions. If one is lucky the linear
programming bounds using truncated regions will still be strong enough to
give the desired inequality. The use of truncation greatly reduces the
number of configurations that must be examined.
Truncation bounds were also used in 
proving Theorem~\ref{th52} above.

We add the following remarks about the construction of
the linear programming problems.
\begin{itemize}
\item[(1)]
For each face $F$ of a given graph type $\sG$ Hales and Ferguson 
construct  a large number of linear
programming constraints in terms of the edge lengths, 
dihedral angles and solid angles
of the polyhedral pieces making up the cluster of $\sD_{HF} (\bv )$ 
over face
$F$ of the graph $\sG$.
The edge lengths, dihedral angles and solid angles are
 variables 
in the linear program. Some of the constraints embody geometric
restrictions that a polyhedron of the given type must satisfy. Others
of them are inequalities relating the weight function of the polyhedron, which
is also a variable in the linear program, to the geometric quantities.
The inequalities bound the score function on the cluster 
(either as a $V$-cell
or as $D$-sets) in terms of these variables. 
There are also some global constraints in the linear program,
 for example that the solid 
angles of the faces around $\bv$ add up to $4 \pi$.

\item[(2)]
The weight function for $D$-sets does not
permit subdivision of the simplex, but the weight function on
the $V$-cell is additive under subdivision, so one can cut such
regions up into smaller pieces if necessary, to get improved
linear programming bounds, by including more stringent constraints.

\item[(3)]
In the linear programming relaxation, a feasible
solution to the constraints need not correspond to any geometrically
constructible vertex $D$-star. All that is required is that
every vertex $D$-star of the particular configuration type correspond
to some feasible point of the linear program.
\end{itemize}

In this fashion one obtains a long list of linear programs, one for
each configuration type, and to 
rule out a map type $\sG$ one needs an upper bound for the 
linear program's objective function strictly below 8 pt.
To rigorously obtain such an upper bound, it suffices to find a feasible 
solution to the dual linear program, and to obtain a good upper
bound one wants the dual feasible solution close to a dual
optimal solution.
The value of the dual $LP$'s objective function is then a 
certified upper bound to the primal $LP$.
To obtain such a certification, it is useful to formulate
the linear programs so that 
the dual linear program has only inequality constraints, with
no equality constraints, so that the feasible region for it is
full-dimensional.
This way, one can guarantee that the dual feasible solution is
{\em strictly} inside the dual feasible region which
facilitates checking feasibility.
This is necessary because the linear program put on the 
computer is only an approximation to the true linear program.
For example, certain constraints of the true $LP$ involve transcendental
numbers like $\pi$, and one considers an approximation.
The effect of these errors is to perturb the {\em objective function} 
of the dual linear program.
Thus a rigorous bound on the effect of these perturbations on the 
upper bound can be obtained in terms of the dual feasible solution.
In this way one can (in principle) get a certified upper bound
\footnote{The Hales proof in the preprints
used a linear programming package CPLEX 
that does not supply such certificates.  Therefore the linear 
programming part of the Ferguson-Hales proof needs to be re-done 
to obtain {\em guaranteed} certificates.}
on the score for a map type $\sG$, using a computer.

The linear programming bounds in the Ferguson-Hales approach above suffice 
to eliminate all map types $\sG$ not ruled out by Theorem~\ref{th52}
except for about 100 ``bad'' cases. These are then handled
by ad hoc methods.
[I am not sure
of the details about how these remaining ``bad'' cases are handled. 
Presumably
they are split into smaller pieces, extra inequalities are generated
somehow, and perhaps specific information on the location of 
vertices more than $\frac{251}{100}$ is incorporated into the linear
programs.]

\section{Concluding Remarks}
\hsp
The Kepler conjecture  appears to be an extraordinarily difficult 
nonlinear optimization problem.
The ``configuration space'' to be optimized over has an extremely 
complicated
structure, of high dimensionality, and the function being optimized is
highly nonlinear and nonconvex, and lacks good monotonicity properties.
The crux of the Hales approach is to select a formulation of an 
optimization problem that can be carried out (mostly by computer) 
in a reasonable length of time.
This led to the Hales-Ferguson choice of an very complicated partition and 
score function, giving an inelegant local inequality, which
however has good decomposition properties in terms of the nonlinear
program. Much of the work in the proof lies in the reductions to
reasonable sized cases, and the use of linear programming
relaxations.
The elimination of the most complicated cases, in Theorem \ref{th52} was
a major accomplishment of this approach. 
The use of Delaunay simplices
to cover most of the volume where density 
is high seems important to the proof and to the choice of score
functions, since simple analytic formulae are available
for tetrahedra.
The Hales - Ferguson proof, assumed correct, is a tour de force 
of nonlinear optimization.

In contrast, the Hsiang approach formulates a relatively elegant local 
inequality, involving only Voronoi domains and a fairly simple 
weight function:
only nearest neighbor regions are counted.
It is conceivable that a rigorous proof of the Hsiang inequality can be 
established, but it very likely will require an enormous 
computer-aided proof of a 
sort very similar to the Hales approach.
Voronoi domains do not seem well suited to computer proof:
they may have 40 or more faces each, and the Hsiang approach requires 
considering up to twenty of them at a time.
A computer-aided proof would likely have to dissect the Voronoi domains 
into pieces, further increasing the size of the problem.

\paragraph{Acknowledgments.} I am indebted to T. Hales for critical
readings of a preliminary version, with many suggestions and corrections.
G. Ziegler provided comments and corrections.

\clearpage
\subsection*{Appendix A. Hales Score Function Formulas}
\hsp
These definitions are taken from in Hales \cite[Section 8]{I} and
Ferguson and Hales \cite[p. 8--11]{FH}.
A tetrahedron $T( l_1, l_2 , l_3 , \ldots, l_6 )$ is
uniquely determined by its six edge lengths $l_i$.
Let the vertices of $T$ be $\bv_0, \bv_1, \bv_2 , \bv_3$ and number
the edges as
\beql{A1}
l_i = \| \bv - \bv_i \|
\quad\mbox{for}\quad
1 \le i \le 3, ~
l_4 = \| \bv_2 - \bv_3 \| ,
l_5 = \| \bv_1 - \bv_3 \|\quad\mbox{and}\quad
l_6 = \| \bv_1 - \bv_2\| ~.
\eeq
We take $\bv_0 = \bzero$ for convenience.

Suppose that the circumcenter $\bw_c = \bw_c (T)$ of $T$ is contained in 
the pointed cone over vertex $\bv_0$, determined by $T$.
Let $\hat{T}_0$ denote the part of the Voronoi cell of $\bv_0$
with respect to the set $\Omega = \{\bv_0, \bv_1, \bv_2 \bv_3 \}$ of 
vertices of $T$ that lies in $T$.
Suppose in addition that the three faces of $T$ containing $\bv_0$
are each {\em non-obtuse} triangles. Then
the set $\hat{T}_0$ subdivides into six pieces, called Rogers
simplices by Hales \cite[p. 31]{I}.
A {\em Rogers simplex} in $T$ is the convex hull of $\bv_0$, the 
midpoint
of an edge emanating from $\bv_0$, the circumcenter of one face of $T$
containing that edge, and the circumcenter $\bw_c = \bw_c (T)$.
If $a$ denotes the half-length of an edge, $b$ the circumradius of a face and
$c= \| \by_c \|$ is the circumradius of $T$ then the associated 
Rogers simplex has shape
\beql{A2}
R(a,b,c) := T(a,b,c, (c^2 -b^1 )^{1/2}, (c^2 -a^2)^{1/2}, (b^2 - a^2 )^{1/2} ) ~,
\eeq
with the positive square root taken.
The intersection of a unit sphere centered at $\bzero$ with 
$R(a,b,c)$ has volume $\frac{1}{3} Sol (\bv_0; R(a,b,c))$, 
where $Sol (\bv_0 ; R(a,b,c))$ denotes the solid angle of 
$R(a,b,c)$ at $\bv_0$,
normalized so that a total solid angle is $4 \pi$.
Set
\beql{A3}
x_i = l_i^2
\eeq
are the squares of the edge lengths.
\begin{lemma}\label{lA1}
The solid angle $Sol ( \bv_0 , T)$ of a tetrahedron 
$T(l_1, l_2, l_3, l_4 , l_5, l_6 )$ is given by
\beql{A4}
Sol ( \bv_0 , T) := 2 ~{\mbox{arccot}} \left( \frac{2A}{\Delta^{1/2}} \right)
\eeq
in which the positive square root of $\Delta$ is taken, 
the value of arccot lies in $[0, \pi ]$, and
\beql{A5}
A(l_1, l_2, l_3, l_4, l_5, l_6) :=
l_1 l_2 l_3 + \frac{1}{2} l_1 (l_2^2 + l_3^2 - l_4^2 ) +
\frac{1}{2} l_2 (l_1^2 + l_3^2 - l_5^2 ) + \frac{1}{2}
l_3 (l_1^2 + l_3^2 - l_6^2 )
\eeq
and
\begin{eqnarray}\label{A6}
\Delta
(l_1, l_2, l_3, l_4, l_5, l_6) & := &
l_1^2 l_4^2 (-l_1^2 + l_2^2 + l_3^2 - l_4^2 + l_5^2 - l_6^2 ) \nonumber \\
&&+ l_2^2 l_5^2 (l_1^2 - l_2^2 + l_3^2 + l_4^2 - l_5^2 + l_6^2 ) \nonumber \\
&&+ l_3^2 l_6^2 (l_2^2 + l_2^2 - l_3^2 + l_4^2 + l_5^2 - l_6^2 ) \nonumber \\
&&- l_2^2 l_3^2 l_4^2 - l_1^2 l_3^2 l_5^2 -l_1^2 l_2^2 l_6^2 - l_4^2 l_5^2 l_6^2 ~.
\end{eqnarray}
\end{lemma}

\paragraph{Definition A.1.}
(i) for a tetrahedron $T(l_1, l_2, \ldots, l_6)$ with vertex $\bv_0$, 
if the circumcenter $\bw_c$ of $T$
falls inside the cone determined by $T$ at $\bv_0$, then we set
\beql{A7}
vor (T, \bv_0) : = 4 \sum_{i=1}^6 \left\{
vol (R_i (a,b,c)) (-\delta_{oct} ) + \frac{1}{3} Sol (R_i, \bv_0)\right\}
\eeq
with
\beql{A8}
vol (R(a,b,c)) := \frac{a(b^2 - a^2)^{1/2} (c^2 - b^2 )^{1/2}}{6} ,
\quad\mbox{for}\quad 1 \le a \le b \le c ~.
\eeq
This formula satisfies $vor (T, \bv_0 ) = 4 \Gamma ( \hat{T}_0)$.

(ii) The six tetrahedra $R_i (a,b,c)$ are still defined even when 
the circumcenter
$\bw_c$ falls outside the cone of $T$ at vertex $\bv_0$, and we still 
take the formula \eqn{A7} to define $vor (T, \bv_0)$, except that both
$vol (R_i (a,b,c))$ and $Sol(R_i,\bv_0)$ are
counted with a negative sign:
each tetrahedron $R_i (a,b,c)$ falls outside $T$, and has no interior 
in common with it.

Hales calls the definition (ii) the ``analytic continuation'' of case (i).
It has a geometric interpretation.

The truncated Voronoi function $vor (T, \bv_0; t)$ of a tetrahedron 
$T$ at vertex $\bv_0$ is intended to measure the compression 
$\Gamma ( \hat{T}_0 \cap B (\bv_0 ; t))$.
Here we have truncated the region $\hat{T}_0$ by removing from it
all points at distance greater than $t$ from $\bv_0$.
We set
\beql{A9}
vor_0 (T, \bv_0) := vor \left(T, \bv_0, \frac{251}{200} \right) ~.
\eeq
The definition
\beql{A10}
vor (T, \bv_0 ; t) := \Gamma ( \hat{T}_0 \cap B (\bv_0 ; t))
\eeq
is valid only when the circumcenter $\bw_c$ of $T$ lies in the 
cone generated from
$T$ at vertex $\bv_0$.
In the remaining case one must construct an analytic representation
analogous to \eqn{A7} for $vor (T, \bv_0; t)$.
This is done in \cite[pages 9--10]{FH}.

\clearpage
\subsection*{Appendix B. References to the Hales Program Results} 
\hsp
This paper was written to state the
Hales-Ferguson local inequality in as simple a way as I could find, 
and does not match 
the order in which things are done in the preprints of Hales and Ferguson. 
Also, the lemmas and theorems stated here are not
all stated in the Hales and Ferguson preprints; some of
them are based on the talks that Hales gave at IAS in January 1999. 
The pointers below indicate where to look in the preprints for 
the results I formulate as lemmas and theorems. 
{\em Warning}:
The Hales-Ferguson partition and scoring function given in \cite{FH},
which are the ones actually used for the proof of the Kepler conjecture, 
differ from those used earlier by Hales in \cite{I} and \cite{II}.

\begin{itemize}
\item[(0)] The idea of considering local inequalities that weight
total area  and covered area by spheres in a ratio $\frac {B}{A}$ that is
not equal to the optimal density occurs in Hales' original approach
based on Delaunay triangulations, see \cite{H1} \cite{H2}.
It also appears in  Hales \cite[Lemma 2.1]{I} 
and in Ferguson and Hales \cite[Proposition 3.14]{FH}. I have 
inserted the parameters $A$ and $B$ in order to include the density
inequality of Hsiang\cite{Hs} in the same framework.
\item[(1)]
{\em Definitions \ref{de41} and \ref{de42}} appear in \cite[page 2]{FH}.
\item[(2)]
{\em Lemma \ref{le41}} is Lemma 1.2 of
\cite{FH}, proved in Lemma 3.5 of \cite{I}.
(The fact that no  vertex of $\Omega$ occurs inside a face of 
a $QL$-tetrahedron or  a $QR$-tetrahedron 
requires additional argument.)
\item[(3)]
{\em Lemma \ref{le42}} follows from Lemma 1.3 of \cite{FH}.
\item[(4)]
{\em Lemma \ref{le43}} is proved in \cite[p.3 bottom]{FH}.

\item[(5)]
{\em Lemma \ref{Nle43}}
is covered in the discussion on \cite[pages 5--6]{FH}.
\item[(6)]
{\em Lemma \ref{Nle44}} (i) -(ii1) are covered in the discussion 
on \cite[page 5]{FH}, including Lemma 1.8 of that paper.
\item[(7)]
The notion of ``tip''  is discussed at length in section 2 of Hales~\cite{II}.
In part II  ``tips'' are not actually reassigned- although this is
mentioned - their existence affects the scoring rule used for the associated
Delaunay simplex which the ``tip'' is associated to.
The rules for moving ``tips'' around to make $V$-cells in the 
Hales-Ferguson approach are
discussed  on \cite[page 8]{FH}. Warning: the way that ``tips'' are
handled in part II and in \cite{FH} may not be the same:
\cite{FH} takes priority. 
\item[(8)]
{\em Lemma \ref{Nle45}} (i) is \cite[Lemma 2.2]{II} and \cite[Lemma 4.17]{FH}.
Facts related to (ii) are discussed in \cite[Sect. 8.6.7]{I}. (For the
second part I do not have a reference.)
(iii) Hales mentioned this in IAS lectures, and sent me a proof sketch, 
which I expanded into the following: Let $S$ be a simplex in the
$D$-system that overlaps a ``tip'' protuding from $\bv$. Say that the ``tip'' 
overlaps by pointing in to $S$ along a face $F$ of $S$. 
Thus  $F$ is a negatively oriented face of $S' =(F, \bv)$,  
which means that the simplex $S'$ is a 
$QR$-tetrahedron or else a $QL$-tetrahedron with spine on $F$.
Suppose first that $S'$ is a $QL$-tetrahedron.
It now follows that $S$ must be
$QL$-tetrahedron with its spine on $F$
by \cite[Lemma 2.2]{FH}.
So $S'$ and $S$ are adjacent $QL$-tetrahedra with spines
on their common face $F$. Now $S$ is in the $D$-system
since $S'$ is in the $D$-system.
Thus the distance of $\bv$ to the vertices in $F$ is at most
$\frac {251}{100},$ since $\bv$ is not on the spine.
We now suppose that  the ``tip'' is not
entirely contained in $S$, and derive a contradiction.
If it isn't contained in $S$, then it crosses out through a face 
$F'$ of $S$.
By the same argument, the distance from $\bv$ to the vertices of
$F'$ is at most $\frac {251}{100}.$ Thus $\bv$ has distances at
most $\frac {251}{100}$ from all vertices of $S$, which is impossible
by \cite[Lemma 1.2]{FH} and \cite[Lemma 1.3]{FH}.
Suppose secondly that  $S'$ is a $QR$-tetrahedron.
Then one shows that $S$ is
also a $QR$-tetrahedron, hence is in the $D$-system. The rest of
the argument goes as before, to the same contradiction.

\item[(9)]
{\em Theorem \ref{th51}}. The main theorem is first
stated as Conjecture 3.15 in \cite[p. 13]{FH}.
It is the theorem asserted to be  proved in \cite{KC}.
\item [(10)]
{\em Lemma \ref{le51}}
(i) appears as \cite[Lemma 3.13]{FH}.
(ii) is a special case of \cite[Lemma 3.13]{FH}
for a quad cluster, which can consist of four congruent $QL$-tetrahedra.
\item[(11)]
The {\em standard regions}~ corresponding to the 
graph $\sG ( \bv )$ are defined on \cite[p. 4]{FH}.
(``Planar map that breaks unit sphere into regions.'')
\item[(12)]
{\em Lemma \ref{le52}}
follows from Lemma 1.6 of \cite{FH}, which implies that crossing 
lines come from $QL$-tetrahedra only.
\item[(13)]
{\em Lemma \ref{le53}} is an immediate consequence of Lemma~\ref{le52}
and \cite[Lemma 1.3]{FH}.
\item[(14)]
{\em Lemma \ref{le54}} appears as \cite[Lemma 3.13]{FH}.
\item[(15)]
{\em Theorem \ref{th52}}
follows from the Corollary to Theorem 4.4 of \cite{IV}.
See also Proposition 7.1 of \cite{III}.
\item[(16)]
{\em Lemma \ref{le55}} and {\em Lemma \ref{le56}}.
These results are briefly stated at the bottom of p. 8 of \cite{FH}.
There are also some relevant details in Hales \cite[Sect. 2.2]{II}.
(I do not know an exact reference for detailed proof.)
\end{itemize}

\begin{center}
\begin{tabular}{rll}
~ & $\underline{\mbox{Hales-Ferguson terminology}}$ & $\underline{\mbox{Terminology in this paper}}$ \\ [+.2in]
(1) & decomposition star ~~~~~~~~~~~~ & vertex D-star \\ [+.2in]
(2) & quasiregular tetrahedron & $QR$-tetrahedron \\ [+.2in]
(3) & quarter & $QL$-tetrahedron \\[+.2in]
(4) & diagonal (of quarter)~~~~~ & spine (of QL-tetrahedron)\\ [+.2in]
(5) & $Q$-system ~~~~~~~~~~~~~& $D$-system  \\ [+.2in]
(6) & score $\sigma(R, \bv)$ ~~~& weight function $\sigma(R, \bv)$ \\ [+.2in]
(7) & standard cluster~~~~~~& cluster
\end{tabular}
\end{center}

\vspace*{.2\baselineskip}
\noindent AT\&T Labs - Research \\
Florham Park, NJ 07932-0971 \\
email: {\tt jcl@research.att.com} \\

\end{document}